\documentclass[12pt,a4paper]{amsart}
\usepackage{amssymb,amsfonts,amsthm,amsmath}

\theoremstyle{plain}
\newtheorem{theorem}{Theorem}

\newtheorem{corollary}{Corollary}

\numberwithin{equation}{section} \numberwithin{theorem}{section}
\numberwithin{lemma}{section} \numberwithin{definition}{section}
\numberwithin{corollary}{section}
\numberwithin{proposition}{section} \textheight =24cm
\begin{document}
\title[Analogues of $q$-series]{The $q$-Dixon sum Dirichlet series analogue}
\author{Geoffrey B Campbell}
\address{Mathematical Sciences Institute,
         The Australian National University,
         Canberra, ACT, 0200, Australia}

\email{Geoffrey.Campbell@anu.edu.au} 

\begin{abstract}
A few years ago, the concept of a $D$-analogue was introduced as a Dirichlet series analogue for the already known
and well researched hypergeometric $q$-series. The $D$-analogue of the $q$-Dixon sum is given here, in the context of seeing a direct comparison with the old Dixon theorem and the $q$ analogue.
\end{abstract}

\maketitle

\section{Introduction} \label{S:intro}

Over more than 120 years' of literature, Dixon's identity (or Dixon's theorem or Dixon's formula) is any of several different but closely related identities proved by A. C. Dixon, some involving finite sums of products of three binomial coefficients, and some evaluating a hypergeometric sum. These identities famously follow from the MacMahon Master theorem, and can now be routinely proved by computer algorithms (Ekhad 1990).
The original identity, from (Dixon 1891), is

  \begin{equation} \label{E:1.0}
 \sum_{k=-a}^{a}(-1)^k \left(\begin{array}{c}
                               2a \\
                               k+a
                             \end{array}\right)^3 = \frac{(3a)!}{(a!)^3}
 ,
  \end{equation}

A generalization, sometimes called Dixon's identity is

  \begin{equation} \label{E:1.0a}
 \sum_{k=-a}^{a}(-1)^k \left(\begin{array}{c}
                               a+b \\
                               a+k
                             \end{array}\right)
                             \left(\begin{array}{c}
                               b+c \\
                               b+k
                             \end{array}\right)
                             \left(\begin{array}{c}
                               c+a \\
                               c+k
                             \end{array}\right) =
                             \frac{(a+b+c)!}{a!b!c!}
 ,
  \end{equation}

 where $a$, $b$, and $c$ are non-negative integers (Wilf 1994, p. 156). The sum on the left can be written as the terminating well-poised hypergeometric series at 1, from (Dixon 1902):

 \begin{equation} \label{E:1.0b}
\left(\begin{array}{c}
                               b+c \\
                               b-a
                             \end{array}\right)
                             \left(\begin{array}{c}
                               c+a \\
                               c-a
                             \end{array}\right)
                                {_{3}F_{2}}(-2a,-a-b,-a-c;1+b-a,1+c-a;1)
  \end{equation}

  and the identity follows as a limiting case (as $a$ tends to an integer) of Dixon's theorem evaluating a well-poised ${_{3}F_{2}}$ generalized hypergeometric series at 1, from (Dixon 1902):

  \begin{equation} \label{E:1.0c}
  _{3}F_{2}(a,b,c;1+a-b,1+a-c;1)=\frac{\Gamma(1+a/2)\Gamma(1+a/2 -b-c)\Gamma(1+a-b)\Gamma(1+a-c )}{\Gamma(1+a)\Gamma(1+a-b-c)\Gamma(1+a/2-b)\Gamma(1+a/2-c)}.\end{equation}

This holds for $\Re(1+ 1/2 -b-c)>0$. As $c$ tends to $-\infty$ it reduces to Kummer's formula for the hypergeometric function $_2F_1$ at $-1$. Dixon's theorem can be deduced from the evaluation of the Selberg integral.

A $q$-analogue of Dixon's formula for the basic hypergeometric series in terms of the $q$-Pochhammer symbol is given by

\begin{equation}
_{4}\Phi_{3}\left[
\begin{matrix}a, -qa^{1/2}, b, c;    & q,{qa^{1/2}}/{bc} \\
               -a^{1/2}, {aq}/b, {aq}/c, &
\end{matrix}
\right]
=\frac{\left(aq,aq/bc,qa^{1/2}/b,qa^{1/2}/c;q\right)}{\left(aq/b,qa^{1/2},qa^{1/2}/{bc};q\right)}, \quad   \label{E:3.4b}
\end{equation}

where $\left| qa^{1/2}/{bc}  \right|<1.$

\section{The $D$ analogue of Dixon's theorem} \label{S:main}

The application of the transform in my paper Campbell \cite{gC2006} to this yields the
\begin{theorem} The Dirichlet series analogue of Dixon's theorem is
\begin{subequations}
\begin{gather}
_{4}\Theta_{3}\left[
\begin{matrix}2a,-\backslash 1+a,b,c;     \gamma,(1+a-b-c)\gamma \\
              -\backslash a,2a+1-b,2a+1-c,
\end{matrix}
\right] \label{E:1.6a}\\
=\frac{\zeta(2a+1-b, 2a+1-c,1+a,1+a-b-c;\gamma)_\infty}
      {\zeta(2a+1,1-b+a,1-c+a,2a+1-b-c;\gamma)_\infty}, \quad and  \label{E:1.6b}
\end{gather}
\end{subequations}
\begin{subequations}
\begin{gather}
_{4}\Theta_{3}\left[
\begin{matrix}m|   \quad   2a,-\backslash 1+a,b,c;    \gamma,(1+a-b-c)\gamma \\
               \quad -\backslash a,2a+1-b,2a+1-c,
\end{matrix}
\right] \label{E:1.7a}\\
=\frac      {J(m| \quad 2a+1,1-b+a,1-c+a,2a+1-b-c;\gamma)_\infty}
            {J(m| \quad 2a+1-b, 2a+1-c,1+a,1+a-b-c;\gamma)_\infty}. \quad   \label{E:1.7b}
\end{gather}
\end{subequations}

\end{theorem}

To make sense of this we utilize the following four definitions from the 2006 paper.

If $\zeta(a)$ is the Riemann zeta function and $\Re \gamma$
chosen such that the functions all exist, define for positive
integers $n$,
\begin{equation}\label{E:2.2}
\zeta(a;\gamma)_{n}=\prod_{k=0}^{n-1}{\zeta\left((a+k)\gamma)\right)}
=\prod_{p}\frac{1}{(p^{-a\gamma};p^{-\gamma})_{n}},
\end{equation}
\begin{equation}\label{E:2.3}
\zeta(a_{1},a_{2},\ldots,a_{r};\gamma)_{n}=\zeta(a_{1};\gamma)_{n}\zeta(a_{2};\gamma)_{n}\cdots\zeta(a_{r};\gamma)_{n}.
\end{equation}

If $\sigma_k(a)$ is the sum of $k$th powers of the divisors of positive integer~$a$  as in (\ref{E:1.6}) then for positive
integers $n$,
\begin{equation}\label{E:2.4}
\sigma_{-\gamma}(a;k)=\prod_{j=0}^{a-2}\frac{\sigma_{-\gamma}(k\prod_{p|k}p^{j})}{\sigma_{-\gamma}(\prod_{p|k}p^{j})},
\quad \text{defined as } 1 \text{ at } a=1, \text{ and}
\end{equation}
\begin{equation}\label{E:2.5}
\sigma_{-\gamma}(a_{1},a_{2},\ldots,a_{r};k)=\sigma_{-\gamma}(a_{1};k)\sigma_{-\gamma}(a_{2};k)\cdots\sigma_{-\gamma}(a_{r};k).
\end{equation}

Furthermore, to cater for the finite products over the primes we define $S_m$ as the set of positive integers that are have the same prime factors as $m$. From these definitions we get particular cases of our above theorem as follows.
\begin{corollary} Some cases of the Dirichlet series analogue of Dixon's sum are:
\begin{gather}
\sum_{k=1}^{\infty}{\frac{\sigma_{-\gamma}(2a,1,1;k) \quad \sigma_{-2\gamma}(2;k)\quad \sigma_{-\gamma}(a;k)}
{\sigma_{-\gamma}(2a,2a,;k) \quad \sigma_{-2\gamma}(a;k)\quad \sigma_{-\gamma}(2;k)\quad k^{(a-1)\gamma}}}  \\
=\frac{\zeta(2a,2a,a-1,a-1;\gamma)_\infty}{\zeta(2a+1,a,a,2a-1;\gamma)_\infty} 
=\frac{\zeta(2a\gamma)\zeta((a-1)\gamma)}{\zeta(a\gamma)\zeta((2a-1)\gamma)}, \notag   \\
\sum_{k\epsilon S_m}{\frac{\sigma_{-\gamma}(2a,1,1;k) \quad \sigma_{-2\gamma}(2;k)\quad \sigma_{-\gamma}(a;k)}
{\sigma_{-\gamma}(2a,2a,;k) \quad \sigma_{-2\gamma}(a;k)\quad \sigma_{-\gamma}(2;k)\quad k^{(a-1)\gamma}}}  \\
=\frac{\sigma_{-\gamma}(\prod_{k|m}{p^{(a-1)}}) \quad \sigma_{-\gamma}(\prod_{k|m}{p^{(2a-2)}})}{\sigma_{-\gamma}(\prod_{k|m}{p^{(2a-1)}}) \quad \sigma_{-\gamma}(\prod_{k|m}{p^{(a-2}})}, \label{E:2.7} \notag  \\
\sum_{k=1}^{\infty}{\frac{\sigma_{-\gamma}(2a,2,2;k) \quad \sigma_{-2\gamma}(1+a;k)\quad \sigma_{-\gamma}(a;k)}
{\sigma_{-\gamma}(2a-1,2a-1,;k) \quad \sigma_{-2\gamma}(a;k) \quad \sigma_{-\gamma}(2;k)\quad k^{(a-3)\gamma}}}   \\
=\frac{\zeta(2a-1,2a-1,a+1,a-3;\gamma)_\infty}{\zeta(2a+1,a,a,2a-1;\gamma)_\infty} \notag \\
=\frac{\zeta((2a-1)\gamma)\zeta((2a)\gamma) \zeta((a-3)\gamma)\zeta((a-2)\gamma)}{\zeta((2a-3)\gamma)\zeta((2a-2)\gamma) \zeta((a-1)\gamma)\zeta((a\gamma)}. \notag \\
\sum_{k\epsilon S_m}{\frac{\sigma_{-\gamma}(2a,2,2;k) \quad \sigma_{-2\gamma}(1+a;k)\quad \sigma_{-\gamma}(a;k)}
{\sigma_{-\gamma}(2a-1,2a-1,;k) \quad \sigma_{-2\gamma}(a;k) \quad \sigma_{-\gamma}(2;k)\quad k^{(a-3)\gamma}}}   \\
=\frac{\sigma_{-\gamma}(\prod_{k|m}{p^{(2a-4)}}) \quad \sigma_{-\gamma}(\prod_{k|m}{p^{(2a-3)}}) \quad \sigma_{-\gamma}(\prod_{k|m}{p^{(a-2)}}) \quad \sigma_{-\gamma}(\prod_{k|m}{p^{(a-1)}})}{\sigma_{-\gamma}(\prod_{k|m}{p^{(2a-2)}}) \quad \sigma_{-\gamma}(\prod_{k|m}{p^{(2a-1}}) \quad \sigma_{-\gamma}(\prod_{k|m}{p^{(a-4)}}) \quad \sigma_{-\gamma}(\prod_{k|m}{p^{(a-3}})}, \label{E:2.7} \notag
 \end{gather}
 \end{corollary}

\end{document}